\pdfoutput=1
\documentclass[english]{article}
\usepackage{geometry}
\usepackage{amsmath}
\usepackage{amsthm}
\usepackage{amssymb}
\usepackage{babel}
\usepackage{authblk}
\begin{document}
\date{}
\title{Technical Report: Infinite Horizon Discrete-Time Linear Quadratic Gaussian Tracking Control Derivation}

\author{Kasra Yazdani}
\author{Matthew Hale}
\affil{Department of  Mechanical and Aerospace Engineering at the University of Florida, Gainesville, FL USA.}
\affil{\texttt {\{kasra.yazdani,matthewhale\}@ufl.edu}}

\maketitle

\section{Introduction\label{sec:Introduction}}

$ $ 

This technical report is meant to accompany \cite{Yazdani2018}, which
relies on the optimal controller for an infinite-horizon LQG tracking
problem in discrete time. Although this controller can be derived
with standard methods, we were unable to find it in the existing literature,
and we therefore derive it here. The notation and the general approach
taken in this report are adapted from the work in \cite[Chapters 1, 4, 5]{Bertsekas2005}.
We adapt that work to derive the optimal controller, which is not
explicitly derived in the form we require it in that reference. Below, we consider stochastic
discrete-time systems of the form

\begin{equation}
x_{k+1}=f_{k}\left(x_{k},u_{k},w_{k}\right),\;k=0,1,\dots,N-1,\label{eq:basic problem}
\end{equation}
where $N$ is the number of timesteps in the control problem. The
state $x_{k}$ is an element of a space $S_{k}$, the control $u_{k}$
is an element of a space $C_{k}$, and the random disturbance $w_{k}$
is an element of a space $D_{k}$. The control $u_{k}$ is constrained
to take values in a given nonempty subset $U\left(x_{k}\right)\subset C_{k}$
which depends on the current state $x_{k}$; that is, $u_{k}\in U_{k}\left(x_{k}\right)$
for all $x_{k}\in S_{k}$ and $k$. We consider the class of policies
that consist of a sequence of functions 
\[
\pi=\left\{ \mu_{0},\dots,\mu_{N-1}\right\} ,
\]
where $\mu_{k}$ maps states $x_{k}$ into controls via $u_{k}=\mu_{k}\left(x_{k}\right),$
and is defined such that $\mu_{k}\left(x_{k}\right)\in U_{k}\left(x_{k}\right)$
for all $x_{k}\in S_{k}$.

The total cost to be minimized is

\[
E\left\{ h_{N}\left(x_{N}\right)+\sum_{k=1}^{N-1}h_{k}\left(x_{k},u_{k},w_{k}\right)\right\} ,
\]
in which $h_{N}\left(x_{N}\right)$ is the terminal cost and $h_{k}\left(x_{k},u_{k},w_{k}\right)$
is the cost incurred at time step $k$. The expected value is with
respect to the joint distribution of the random variables involved,
which will be detailed below. 

The main technique used in this report is dynamic programming which
is based on the the \emph{principle of optimality }introduced by Bellman
\cite{Bellman2003}, which we will briefly introduce now; see \cite[Page 18]{Bertsekas2005}
for a more complete exposition.

\emph{Principle of Optimality: }Let $\pi^{\ast}=\left\{ \mu_{0}^{\ast},\mu_{1}^{\ast},\dots,\mu_{N-1}^{\ast}\right\} $
be an optimal policy for the basic problem in Equation $\left(\ref{eq:basic problem}\right)$,
and assume that, when using $\pi^{\ast}$, the state $x_{i}$ occurs
at time $i$ with positive probability. Consider the subproblem in
which the state is $x_{i}$ at time $i$ and we wish to minimize the
``cost-to-go'' from time $i$ to time $N$, namely,

\begin{equation}
E\left\{ h_{N}\left(x_{N}\right)+\sum_{k=i}^{N-1}h_{k}\left(x_{k},\mu_{k}\left(x_{k}\right),w_{k}\right)\right\} .
\end{equation}
The principle of optimality states that the truncated policy $\left\{ \mu_{i}^{\ast},\mu_{i+1}^{\ast},\dots,\mu_{N-1}^{\ast}\right\} $
is optimal for this subproblem.

We now formally define the dynamic programming algorithm. 

\emph{DP Algorithm: }For every initial state $x_{0}$, the optimal
cost $J^{*}\left(x_{0}\right)$ of the basic problem is equal to $J_{0}\left(x_{0}\right)$,
given by the last step of the following algorithm, which proceeds
backward in time from period $N-1$ to period $0$:

\begin{equation}
J_{N}\left(x_{N}\right)=h_{N}\left(x_{N}\right),
\end{equation}

\begin{equation}
J_{k}\left(x_{k}\right)=\min_{u_{k}\in U_{k}\left(x_{k}\right)}\underset{w_{k}}{E}\left\{ h_{k}\left(x_{k},u_{k},w_{k}\right)+J_{k+1}\left(f_{k}\left(x_{k},u_{k},w_{k}\right)\right)\right\} ,\quad k=0,1,\dots,N-1,\label{eq:DP Cost Function}
\end{equation}
where the expectation is taken with respect to the probability distribution
of $w_{k}$, which can depend on $x_{k}$ and $u_{k}$. Furthermore,
if $u^{\ast}=\mu^{\ast}\left(x_{k}\right)$ minimizes the right side
of Equation (\ref{eq:DP Cost Function}) for each $x_{k}$ and $k$,
the policy $\pi^{\ast}=\left\{ \mu_{0}^{\ast},\mu_{1}^{\ast},\dots,\mu_{N-1}^{\ast}\right\} $
is optimal.

\section{Basic Problem with Imperfect State Information\label{sec:Imperfect State Information Problem-1}}

$ $

Consider the basic problem introduced in Section \ref{sec:Introduction}
where the controller has access to the noisy observations $z_{k}$
of the form

\begin{align}
z_{k} & =h_{k}\left(x_{k},u_{k-1},v_{k}\right),\quad k=1,2,\dots,N-1\nonumber \\
z_{0} & =h_{0}\left(x_{0},v_{0}\right).
\end{align}

The observation $z_{k}$ belongs to a given observation space $Z_{k}$.
The observation noise $v_{k}$ belongs to a given space $V_{k}$ which
is characterized by a given probability distribution

\[
P_{v_{k}}\left(\cdot\mid x_{k},\dots,x_{0},u_{k-1},\dots,u_{0},w_{k-1},\dots,w_{0},v_{k-1},\dots,v_{0}\right).
\]

The initial state $x_{0}$ is random and is characterized by $P_{x_{0}}$.
The probability distribution of the stochastic disturbances in the
system, which is denoted by $P_{w_{k}}(\cdot\mid x_{k},u_{k})$, is
given, and $w_{k}$ does not depend on prior disturbances $w_{0},\dots,w_{k-1},v_{0},\dots,v_{k-1}$,
but it may depend explicitly on $x_{k}$ and $u_{k}$. The control
$u_{k}$ is constrained to take values from a given nonempty subset
$U_{k}$ of the control space $C_{k}$. We are assuming that this
subset does not depend on $x_{k}$.

The information available to the controller at time $k$ is denoted
$I_{k}$ and is called the information vector, formally defined as

\begin{align}
I_{k} & =(z_{0},z_{1},\dots,z_{k},u_{0},u_{1},\dots,u_{k-1}),\quad k=1,2,\dots,N-1\nonumber \\
I_{0} & =z_{0}.
\end{align}

We consider the class of policies consisting of a sequence of functions
$\pi=\left\{ \mu_{0},\mu_{1},\dots,\mu_{N-1}\right\} ,$ where each
function $\mu_{k}$ maps the information vector $I_{k}$ into the
control space $C_{k},$ and

\[
\mu_{k}\left(I_{k}\right)\in U_{k}\text{ for all }I_{k},\quad k=0,1,\dots,N-1,
\]
i.e., such policies are admissible. We denote the reference trajectories by $\bar{x}_k\in T_k$, enabling us to express desire to keep the states of the system close to a given trajectory. We want to find the policy $\pi=\left\{ \mu_{0},\mu_{1},\dots,\mu_{N-1}\right\}$
that minimizes the cost function 

\[
J_{\pi}=\underset{\underset{k=0,\dots,N-1}{x_{0},w_{k},v_{k}}}{E}\left\{ h_{N}\left(x_{N},\bar{x}_{N}\right)+\sum_{k=0}^{N-1}h_{k}\left(x_{k},\mu_{k}\left(I_{k}\right),\omega_{k},\bar{x}_{k}\right)\right\} 
\]
subject to the system dynamics

\[
x_{k+1}=f_{k}\left(x_{k},\mu_{k}\left(I_{k}\right),w_{k}\right),\quad k=1,2,\dots,N-1
\]
and the measurement equation

\begin{align*}
z_{k} & =h_{k}\left(x_{k},\mu_{k-1}\left(I_{k-1}\right),v_{k}\right),\quad k=1,2,\dots,N-1\\
z_{0} & =h_{0}\left(x_{0},v_{0}\right).
\end{align*}

\subsection{Reformulation as Perfect State Information Problem}

$ $

In this section, we will show the reduction from imperfect (i.e.,
noisy state information) to perfect state information. We define a
system whose state at time $k$ is the set of all variables whose
values can benefit the controller when making the $k^{th}$ decision.
By definition,

\begin{align}
I_{k+1} & =\left(I_{k},z_{k+1},u_{k}\right),\quad k=1,2,\dots,N-2,\nonumber \\
I_{0} & =z_{0}.
\end{align}
These equations can be viewed as a new system where the state of the
system is $I_{k}$, control is $u_{k}$, and the measurement $z_{k+1}$
can be treated as a random disturbance. We then have

\[
P\left(z_{k+1}\mid I_{k},u_{k}\right)=P\left(z_{k+1}\mid I_{k},u_{k},z_{0},z_{1},\dots,z_{k}\right).
\]
Therefore, since we already have $z_{0},z_{1},\dots,z_{k}$ as part
of the information vector, the probability distribution of $z_{k+1}$
depends explicitly on the state $I_{k}$ and control of the new system,
$u_{k}$, and not on the prior disturbances $z_{0},z_{1},\dots,z_{k}$.
In order to reformulate the cost function in terms of the variables
of the new system we write

\[
E\left\{ h_{k}\left(x_{k},u_{k},w_{k}\right)\right\} =E\left\{ \underset{x_{k},w_{k}}{E}\left\{ h_{k}\left(x_{k},u_{k},w_{k}\right)\mid I_{k},u_{k}\right\} \right\} .
\]
Hence, the cost per stage as a function of the new state $I_{k}$
and $u_{k}$ is 

\[
\tilde{h}_{k}\left(I_{k},u_{k}\right)=\underset{x_{k},w_{k}}{E}\left\{ h_{k}\left(x_{k},u_{k},w_{k}\right)\mid I_{k},u_{k}\right\} .
\]
The dynamic programming algorithm for $k=0,1,\dots,N-2$ is then written
as

\begin{equation}
J_{k}\left(I_{k}\right)=\underset{u_{N-1}\in U_{N-1}}{\min}\left[\underset{x_{k},w_{k},z_{k+1}}{E}\left\{ h_{k}\left(x_{k},u_{k},w_{k}\right)+J_{k+1}\left(I_{k},z_{k+1},u_{k}\right)\mid I_{k},u_{k}\right\} \right],\label{eq:DP for the stochastic system}
\end{equation}
where $J_{k+1}\left(I_{k},z_{k+1},u_{k}\right)=J_{k+1}\left(I_{k+1}\right)$.
In Equation $\left(\ref{eq:DP for the stochastic system}\right)$,
it is important to note that the expected value is taken with respect
to the observation $z_{k+1}$ as well, since the information vector
$I_{k+1}$ is explicitly dependent on $z_{k+1}$. The cost function
for stage $k=N-1$ is

\begin{align}
J_{N-1}\left(I_{N-1}\right) & =\underset{u_{N-1}\in U_{N-1}}{\min}\nonumber \\
 & \left[\underset{x_{N-1},w_{N-1}}{E}\left\{ h_{N}\left(f_{N-1}\left(x_{N-1},u_{N-1},w_{N-1}\right)\right)+h_{N-1}\left(x_{N-1},u_{N-1},w_{N-1}\right)\mid I_{N-1},u_{N-1}\right\} \right],
\end{align}
where $h_{N}\left(f_{N-1}\left(x_{N-1},u_{N-1},w_{N-1}\right)\right)=h_{N}\left(x_{N}\right)$.

\subsection{Linear Quadratic Gaussian Tracking Controller Derivation}

$ $

We now apply the DP principle to solve for the linear quadratic cost
function. Consider the linear system equation

\begin{equation}
x_{k+1}=A_{k}x_{k}+B_{k}u_{k}+w_{k},\quad k=1,2,\dots,N-1,
\end{equation}
and the quadratic cost function

\begin{equation}
E\left\{ \left(x_{N}-\bar{x}_{N}\right)^{T}Q_{N}\left(x_{N}-\bar{x}_{N}\right)+\sum_{k=0}^{N-1}\left(\left(x_{k}-\bar{x}_{k}\right)^{T}Q_{k}\left(x_{k}-\bar{x}_{k}\right)+u_{k}^{T}R_{k}u_{k}\right)\right\} ,
\end{equation}
where the reference trajectory $\bar{x}_{k}$ is an element of space
$T_{k}$. At the beginning of period $k$ we get an observation $z_{k}$
of the form

\begin{equation}
z_{k}=C_{k}x_{k}+v_{k},
\end{equation}
where the vectors $v_{k}$ are mutually independent, and independent
from $w_{k}$ and $x_{0}$ as well. The noise terms $w_{k}$ are assumed
to be independent, zero mean, and have a finite variance. $C_{k}$
is a given $s\times n$ matrix. We also assume $Q_{k}\succ0$ and
$R_{k}\succ0$. Then, the cost function for step $N-1$ can be expressed
as 

\begin{align}
J_{N-1}\left(I_{N-1}\right) & =\underset{u_{N-1}}{\min}\Bigg[\underset{x_{N-1},w_{N-1}}{E}\bigg\{\left(x_{N-1}-\bar{x}_{N-1}\right)^{T}Q_{N-1}\left(x_{N-1}-\bar{x}_{N-1}\right)+u_{N-1}^{T}R_{N-1}u_{N-1}\nonumber \\
 & +\left(Ax_{N-1}+Bu_{N-1}+w_{N-1}-\bar{x}_{N}\right)^{T}Q_{N}\left(Ax_{N-1}+Bu_{N-1}+w_{N-1}-\bar{x}_{N}\right)\mid I_{N-1},u_{N-1}\bigg\}\Bigg].
\end{align}

Since we have $E\left\{ w_{N-1}\mid I_{N-1}\right\} =E\left\{ w_{N-1}\right\} =0$,
we can write

\begin{align}
J_{N-1}\left(I_{N-1}\right) & =\underset{x_{N-1}}{E}\Bigg\{\left(x_{N-1}-\bar{x}_{N-1}\right)^{T}Q_{N-1}\left(x_{N-1}-\bar{x}_{N-1}\right)\nonumber \\
 & +\left(Ax_{N-1}\right)^{T}Q_{N}\left(Ax_{N-1}\right)-2\left(Ax_{N-1}\right)^{T}Q_{N}\bar{x}_{N}+\bar{x}_{N}^{T}Q_{N}\bar{x}_{N}\mid I_{N-1}\bigg\}\nonumber \\
 & +\underset{u_{N-1}}{\min}\left[u_{N-1}^{T}\left(B^{T}Q_{N}B+R_{N-1}\right)u_{N-1}+2\left(AE\left\{ x_{N-1}\mid I_{N-1}\right\} -\bar{x}_{N}\right)^{T}Q_{N}Bu_{N-1}\right]\nonumber \\
 & +\underset{w_{N-1}}{E}\left\{ w_{N-1}^{T}Q_{N}w_{N-1}\right\} .\label{eq:Cost at step N-1 before taking min}
\end{align}

In order to minimize the cost for the last period, we set $\frac{dJ_{N-1}\left(I_{N-1}\right)}{du_{N-1}}=0$
and find

\begin{align}
u_{N-1}^{*} & =\mu_{N-1}^{*}\left(I_{N-1}\right)\nonumber \\
 & =-\left(B^{T}Q_{N}B+R_{N-1}\right)^{-1}B^{T}Q_{N}\left(AE\left\{ x_{N-1}\mid I_{N-1}\right\} -\bar{x}_{N}\right).
\end{align}

Substituting this into Equation$\left(\ref{eq:Cost at step N-1 before taking min}\right)$
provides

\begin{align*}
J_{N-1}\left(I_{N-1}\right) & =\underset{x_{N-1}}{E}\Bigg\{\left(x_{N-1}-\bar{x}_{N-1}\right)^{T}Q_{N-1}\left(x_{N-1}-\bar{x}_{N-1}\right)\\
+ & \left(Ax_{N-1}-\bar{x}_{N}\right)^{T}Q_{N}\left(Ax_{N-1}-\bar{x}_{N}\right)\\
+ & \Bigg(\left(B^{T}Q_{N}B+R_{N-1}\right)^{-1}B^{T}Q_{N}\left(AE\left\{ x_{N-1}\mid I_{N-1}\right\} -\bar{x}_{N}\right)\Bigg)^{T}\\
.\; & \left(B^{T}Q_{N}B+R_{N-1}\right)\left(B^{T}Q_{N}B+R_{N-1}\right)^{-1}B^{T}Q_{N}\left(AE\left\{ x_{N-1}\mid I_{N-1}\right\} -\bar{x}_{N}\right)\\
- & 2\left(AE\left\{ x_{N-1}\mid I_{N-1}\right\} -\bar{x}_{N}\right)^{T}Q_{N}B\left(B^{T}Q_{N}B+R_{N-1}\right)^{-1}B^{T}Q_{N}\left(AE\left\{ x_{N-1}\mid I_{N-1}\right\} -\bar{x}_{N}\right)\mid I_{N-1}\Bigg\}\\
+ & \underset{w_{N-1}}{E}\left\{ w_{N-1}^{T}Q_{N}w_{N-1}\right\} .
\end{align*}

Toward completing the square, we add and subtract $\underset{x_{N-1}}{E}\Bigg\{\left(Ax_{N-1}-\bar{x}_{N}\right)^{T}\underbar{P}_{N-1}\left(Ax_{N-1}-\bar{x}_{N}\right)\Bigg\}$
to get

\begin{align*}
J_{N-1}\left(I_{N-1}\right) & =\underset{x_{N-1}}{E}\Bigg\{\left(x_{N-1}-\bar{x}_{N-1}\right)^{T}Q_{N-1}\left(x_{N-1}-\bar{x}_{N-1}\right)+\left(Ax_{N-1}-\bar{x}_{N}\right)^{T}Q_{N}\left(Ax_{N-1}-\bar{x}_{N}\right)\mid I_{N-1}\Bigg\}\\
+ & \underset{x_{N-1}}{E}\Bigg\{\left(AE\left\{ x_{N-1}\mid I_{N-1}\right\} -\bar{x}_{N}\right)^{T}\underbar{P}_{N-1}\left(AE\left\{ x_{N-1}\mid I_{N-1}\right\} -\bar{x}_{N}\right)\\
- & 2\left(Ax_{N-1}-\bar{x}_{N}\right)^{T}\underbar{P}_{N-1}\left(AE\left\{ x_{N-1}\mid I_{N-1}\right\} -\bar{x}_{N}\right)\mid I_{N-1}\Bigg\}\\
+ & \underset{x_{N-1}}{E}\Bigg\{\left(Ax_{N-1}-\bar{x}_{N}\right)^{T}\underbar{P}_{N-1}\left(Ax_{N-1}-\bar{x}_{N}\right)\Bigg\}-\underset{x_{N-1}}{E}\Bigg\{\left(Ax_{N-1}-\bar{x}_{N}\right)^{T}\underbar{P}_{N-1}\left(Ax_{N-1}-\bar{x}_{N}\right)\Bigg\}\\
+ & \underset{w_{N-1}}{E}\left\{ w_{N-1}^{T}Q_{N}w_{N-1}\right\} ,
\end{align*}
where we define $\underbar{P}_{N-1}\triangleq Q_{N}B\left(B^{T}Q_{N}B+R_{N-1}\right)^{-1}B^{T}Q_{N}$.
Completing the square, we find

\begin{align}
J_{N-1}\left(I_{N-1}\right) & =\underset{x_{N-1}}{E}\Bigg\{\left(x_{N-1}-\bar{x}_{N-1}\right)^{T}Q_{N-1}\left(x_{N-1}-\bar{x}_{N-1}\right)\nonumber \\
 & +\left(Ax_{N-1}-\bar{x}_{N}\right)^{T}\left(Q_{N}-\underbar{P}_{N-1}\right)\left(Ax_{N-1}-\bar{x}_{N}\right)\mid I_{N-1}\Bigg\}\nonumber \\
 & +\underset{x_{N-1}}{E}\Bigg\{\left[\left(Ax_{N-1}-\bar{x}_{N}\right)-\left(AE\left\{ x_{N-1}\mid I_{N-1}\right\} -\bar{x}_{N}\right)\right]^{T}\underbar{P}_{N-1}\nonumber \\
 & .\left[\left(Ax_{N-1}-\bar{x}_{N}\right)-\left(AE\left\{ x_{N-1}\mid I_{N-1}\right\} -\bar{x}_{N}\right)\right]\mid I_{N-1}\Bigg\}\nonumber \\
 & +\underset{w_{N-1}}{E}\left\{ w_{N-1}^{T}Q_{N}w_{N-1}\right\} .\label{eq:Cost for step N-1}
\end{align}
Here, the terms containing $\bar{x}_{N}$ cancel each other out in
the $2^{nd}$ and $3^{rd}$ line. Inductively, we write the cost function
for step $N-2$ as

\begin{align}
J_{N-2}\left(I_{N-2}\right) & =\underset{u_{N-2}}{\min}\Bigg[\underset{x_{N-2},w_{N-2},x_{N-1}}{E}\Bigg\{\left(x_{N-2}-\bar{x}_{N-2}\right)^{T}Q_{N-2}\left(x_{N-2}-\bar{x}_{N-2}\right)\nonumber \\
 & +u_{N-2}^{T}R_{N-2}u_{n-2}+J_{N-1}\left(I_{N-1}\right)\mid I_{N-2},u_{N-2}\bigg\}\Bigg],
\end{align}
and substituting for $J_{N-1}$ from Equation (\ref{eq:Cost for step N-1})
gives

\begin{align}
J_{N-2}\left(I_{N-2}\right) & =\underset{x_{N-2}}{E}\Bigg\{\left(x_{N-2}-\bar{x}_{N-2}\right)^{T}Q_{N-2}\left(x_{N-2}-\bar{x}_{N-2}\right)\mid I_{N-2}\Bigg\}\nonumber \\
+ & \underset{u_{N-2}}{\min}\Bigg[u_{N-2}^{T}R_{N-2}u_{n-2}+\underset{x_{N-1}}{E}\Bigg\{ x_{N-1}^{T}\left(Q_{N-1}+A^{T}Q_{N}A-A^{T}\underbar{P}_{N-1}A\right)x_{N-1}\nonumber \\
+ & \bar{x}_{N-1}^{T}Q_{N-1}\bar{x}_{N-1}+\bar{x}_{N}^{T}\left(Q_{N}-\underbar{P}_{N-1}\right)\bar{x}_{N}\nonumber \\
- & 2\bar{x}_{N-1}^{T}Q_{N-1}x_{N-1}-2\bar{x}_{N}^{T}\left(Q_{N}-\underbar{P}_{N-1}\right)Ax_{N-1}\mid I_{N-1}\Bigg\}\Bigg]\nonumber \\
+ & \underset{x_{N-1}}{E}\left\{ \left[A\left(x_{N-1}-E\left\{ x_{N-1}\mid I_{N-1}\right\} \right)\right]^{T}\underbar{P}_{N-1}\left[A\left(x_{N-1}-E\left\{ x_{N-1}\mid I_{N-1}\right\} \right)\right]\mid I_{N-1}\right\} \nonumber \\
+ & \underset{w_{N-1}}{E}\left\{ w_{N-1}^{T}Q_{N}w_{N-1}\right\} .\label{eq: Cost for  N-2}
\end{align}
It is important to note that based on \cite[Lemma 5.2.1]{Bertsekas2005},
the term in line 4 of Equation (\ref{eq: Cost for  N-2}) is not
in the above minimization with respect to $u_{N-2}$. Next, define

\begin{equation}
K_{N-1}\triangleq Q_{N-1}+A^{T}Q_{N}A-P_{N-1}
\end{equation}

\begin{equation}
P_{N-1}\triangleq A^{T}Q_{N}B\left(B^{T}Q_{N}B+R_{N-1}\right)^{-1}B^{T}Q_{N}A,
\end{equation}
and then set $\frac{dJ_{N-2}\left(I_{N-2}\right)}{du_{N-2}}=0$ to
find

\begin{multline*}
2u_{N-2}^{T}R_{N-2}\\
+E\Bigg\{2\left(Ax_{N-2}+Bu_{N-2}+w_{N-2}\right)^{T}K_{N-1}B-2B^{T}Q_{N-1}\bar{x}_{N-1}-2B^{T}A^{T}\left(Q_{N}-\underbar{P}_{N-1}\right)\bar{x}_{N}\mid I_{N-1}\Bigg\}=0,
\end{multline*}
where we used $x_{N-1}=Ax_{N-2}+Bu_{N-2}+w_{N-2}$. Simplifying this
equation results in

\begin{multline*}
u_{N-2}^{T}R_{N-2}\\
+u_{N-2}^{T}B^{T}K_{N-1}B+B^{T}K_{N-1}AE\left\{ x_{N-1}\mid I_{N-1}\right\} -B^{T}Q_{N-1}\bar{x}_{N-1}-B^{T}A^{T}\left(Q_{N}-\underbar{P}_{N-1}\right)\bar{x}_{N}\mid I_{N-1}=0,
\end{multline*}
and solving for $u_{N-2}^{*}$ provides

\[
u_{N-2}^{*}=-\left(R_{N-2}+B^{T}K_{N-1}B\right)^{-1}B^{T}\left(K_{N-1}AE\left\{ x_{N-1}\mid I_{N-1}\right\} -Q\bar{x}_{N-1}-A\left(Q_{N}-\underbar{P}_{N-1}\right)\bar{x}_{N}\right).
\]

Continuing this strategy and putting it into recursive form, we have
the optimal control

\begin{equation}
u_{k}^{*}=-\left(R_{k}+B^{T}K_{k+1}B\right)^{-1}B^{T}\left(K_{k+1}AE\left\{ x_{k}\mid I_{k}\right\} +g_{k+1}\right),
\end{equation}
where

\begin{equation}
g_{k}=A^{T}\left[I-K_{k+1}B\left(R+B^{T}K_{k+1}B\right)^{-1}B^{T}\right]g_{k+1}-Q\bar{x}_{k}
\end{equation}
\begin{equation}
g_{N}=-Q_{N}\bar{x}_{N}
\end{equation}

\begin{equation}
K_{N}=Q_{N}
\end{equation}

\begin{equation}
P_{k}=A^{T}K_{k+1}^{T}B\left(B^{T}K_{k+1}B+R_{k}\right)^{-1}B^{T}K_{k+1}A
\end{equation}

\begin{equation}
K_{k}=A^{T}K_{k+1}A-P_{k}+Q_{k}.
\end{equation}

The authors in \cite{Kwon2005} investigate a closely related problem
where they assume access to non-stochastic state knowledge and absence
of process noise. They derive a closed form solution to the linear
quadratic tracking control problem using the minimum principle of
Pontryagin, and the results in this report are consistent with the
results in \cite{Kwon2005}.
\begin{flushright}
\par\end{flushright}

\section{Quantifying Increase in Cost When Tracking a Stochastic Reference Trajectory}
A common extension to the problem solved in Section
\ref{sec:Imperfect State Information Problem-1} is the case where
the problem is solved in steady state across an infinite horizon.
In particular, we are interested in infinite horizon LQG tracking
control problems in which the reference trajectory is a stochastic
process.

The cost function for finite-time horizon LQG tracking a stochastic
reference trajectory is defined as 

\[
J_N^{LQGSTC}(x,u)=E\left\{ \left(x_{N}-\tilde{x}_{N}\right)Q_{N}\left(x_{N}-\tilde{x}_{N}\right)+\sum_{k=0}^{N-1}\left(\left(x_{k}-\tilde{x}_{k}\right)^{T}Q_{k}\left(x_{k}-\tilde{x}_{k}\right)+\tilde{u}_{k}^{T}R_{k}\tilde{u}_{k}\right)\right\},
\]
where the stochastic reference trajectory $\tilde{x}_{k}$ is defined
as 

\[
\tilde{x}_{k}=\bar{x}_{k}+\bar{w}_{k}.
\]
We further develop the LQG tracking control problem
by analyzing the increase in the cost due to the added noise to the
reference trajectories. We assume that $\bar{x}$ has a limiting value
i.e., 
\[
\lim_{k\rightarrow\infty}\bar{x}_{k}=\bar{x}
\]
is finite. 

In an average cost per stage problem, temporary increases in cost average
out over infinite time horizons, and one can make several simplifying
modifications to infinite-horizon problem statements without changing
the optimal controller. One such modification comes from replacing
each element of the reference trajectory by the limiting value
of the reference trajectory.  
Doing  so, the cost incurred in the LQG tracking control problem over an infinite time horizon when tracking a stochastic reference trajectory is 
equal to

\begin{equation}
J^{IHLQGSTC}(x,u)=\lim_{N\rightarrow\infty}\frac{1}{N}\underset{w_{k},\tilde{x}}{E}\left\{ \sum_{k=0}^{N-1}\left(\left(x_{k}-\tilde{x}\right)^{T}Q\left(x_{k}-\tilde{x}\right)+\tilde{u}_{k}^{T}R\tilde{u}_{k}\right)\right\} ,\label{eq:Infinite Horizon Cost}
\end{equation}
where $\tilde{x}=\bar{x}+\bar{w}$ and $\bar{w}$ is a noise term satisfying the same assumptions as $\bar{w}_k$ above. Because this is an average 	cost-per-stage problem, we can use a time-invariant controller based on only $\tilde{x}$, rather than the entire sequence $\{\tilde{x}(k)\}_{k\in {\mathbb{N}}}$, and any increases in cost average out over time.

The optimal control value  $\tilde{u}^*_k$ when tracking a stochastic reference trajectory is defined as 

\begin{equation}
\tilde{u}_{k}^{*}=-\left(R+B^{T}KB\right)^{-1}B^{T}\left(KAE\left\{ x_{k}\mid I_{k}\right\} +\tilde{g}\right),
\end{equation}
where

\begin{equation}\label{eq:gtilde}
\tilde{g}=A^{T}\left[I-KB\left(R+B^{T}KB\right)^{-1}B^{T}\right]\tilde{g}-Q\tilde{x},
\end{equation}
and where $K$ is the unique positive semi-definite solution to the
algebraic Riccati equation

\begin{equation}
K=A^{T}KA-A^{T}K^{T}B\left(B^{T}KB+R\right)^{-1}B^{T}KA+Q.
\end{equation}

One might wonder, how will tracking a stochastic reference trajectory affect the cost incurred in infinite time relative to tracking a non-stochastic reference trajectory. We start by solving Equation~\eqref{eq:gtilde} via 

\begin{align}
\tilde{g}&=-\left( I - A^{T}\left(I-KB\left(R+B^{T}KB\right)^{-1}B^{T}\right) \right)^{-1}Q \tilde{x}\\
& = -\left( I - A^{T}\left(I-KB\left(R+B^{T}KB\right)^{-1}B^{T}\right) \right)^{-1} Q(\bar{x}+\bar{w})\\
& = g - \left( I - A^{T}\left(I-KB\left(R+B^{T}KB\right)^{-1}B^{T}\right) \right)^{-1} Q\bar{w}.
\end{align}
Substituting $\tilde{g}$ in $\tilde{u}_{k}^{*}$, we find 
\begin{align}
\tilde{u}_{k}^{*}&=-\left(R+B^{T}KB\right)^{-1}B^{T}\left(KAE\left\{ x_{k}\mid I_{k}\right\} +\tilde{g}\right)\\
& = -\left(R+B^{T}KB\right)^{-1}B^{T}\left(KAE\left\{ x_{k}\mid I_{k}\right\}\right) -\left(R+B^{T}KB\right)^{-1}B^{T}\tilde{g}\\
& =  -\left(R+B^{T}KB\right)^{-1}B^{T}\left(KAE\left\{ x_{k}\mid I_{k}\right\}\right) \\
&-\left(R+B^{T}KB\right)^{-1}B^{T}\left[ g -\left( I - A^{T}\left(I-KB\left(R+B^{T}KB\right)^{-1}B^{T}\right) \right)^{-1} Q \bar{w} \right] \\
& = u_{k}^{*} + \left(R+B^{T}KB\right)^{-1}B^{T} \left( I - A^{T}\left(I-KB\left(R+B^{T}KB\right)^{-1}B^{T}\right) \right)^{-1} Q\bar{w},
\end{align}
where here~$u^*_k$ denotes the infinite-horizon analog  to the
controller developed in Section~2. 

Next, define 
\begin{equation}
F = \left(R+B^{T}KB\right)^{-1}B^{T} \left( I - A^{T}\left(I-KB\left(R+B^{T}KB\right)^{-1}B^{T}\right) \right)^{-1}Q.
\end{equation}
The cost at each time $k$ is

\begin{align}
E \left\lbrace \left(x_{k}-\tilde{x}\right)^{T}Q_k\left(x_{k}-\tilde{x}\right)+\tilde{u}_{k}^{T}R_k\tilde{u}_{k} \right\rbrace, 
\end{align}
and substituting for $\tilde{u}_k^{*}$ gives

\begin{align}
E \left\lbrace \left(x_{k}-\tilde{x}\right)^{T}Q_k\left(x_{k}-\tilde{x}\right)+ \left( u^*_k + F \bar{w} \right)^T R_k \left( u^*_k + F \bar{w} \right)  \right\rbrace.
\end{align}
Simplifying shows the above to be equal to 

\begin{align}
\underset{x_k,u_k}{E}\left\{\left(x_{k}-\bar{x}_{k}\right)^{T}Q_{k}\left(x_{k}-\bar{x}_{k}\right)+u_{k}^{T}R_{k}u_{k}\right\} + E\left\{ \bar{w}^{T}Q_k\bar{w}\right\} + E \left\lbrace (F\bar{w})^T R_k (F\bar{w}) \right\rbrace.
\end{align}
Therefore, for $N$ steps we can write

\begin{equation}
J_N^{LQGSTC}(x_0)=J_N^{LQGTC}(x_0)+\sum_{k=0}^{N-1}\bigg( E\left\{ \bar{w}^{T}Q_{k}\bar{w}\right\}+ E \left\lbrace (F\bar{w})^T R_k (F\bar{w}) \right\rbrace \bigg) ,\label{eq:Cost Difff btwn stochastic and deterministic ref}
\end{equation}
where $J_N^{LQGSTC}(x_0)$ denotes the optimal cost of the LQGSTC problem and $J_N^{LQGTC}(x_0)$ denotes the  optimal cost of the LQGTC problem. Both involve $N$ stages, initial state $x_{0}$, and zero terminal cost.  Supposing that both $Q$ and $R$ are constant in time, taking the limit as $N\rightarrow \infty$ gives
\begin{equation}
J^{IHLQGSTC}(x_{0})=\lim_{N\rightarrow\infty}\frac{1}{N}\left(J_N^{LQGTC}(x_{0})+\sum_{k=0}^{N-1}\bigg( E\left\{ \bar{w}^{T}Q\bar{w}\right\}+ E \left\lbrace (F\bar{w})^T R (F\bar{w}) \right\rbrace \bigg) \right),\label{eq:IHLQGSTCbefore cancellation}
\end{equation}
which is equal to
\begin{equation}
J^{IHLQGSTC}(x_{0})=\lim_{N\rightarrow\infty}\left(\frac{1}{N}J_N^{LQGTC}(x_{0})\right)+E\left\{ \bar{w}^{T}Q\bar{w}\right\} + E \left\lbrace (F\bar{w})^T R (F\bar{w}) \right\rbrace   .
\end{equation}

\newpage
\bibliographystyle{plain}
\bibliography{Kasra_bibtex}

\begin{thebibliography}{1}

\bibitem{Bellman2003}
Richard Bellman.
\newblock {\em {Dynamic Programming}}.
\newblock Dover Publications, 2003.

\bibitem{Bertsekas2005}
Dimitri~P. Bertsekas.
\newblock {\em {Dynamic Programming and Optimal Control}}.
\newblock Athena Scientific, 3 edition, 2005.

\bibitem{Kwon2005}
Wook~Hyun Kwon and Soo~Hee Han.
\newblock {\em {Receding Horizon Control: Model Predictive Control for State
  Models}}.
\newblock {Springer-Verlag}, 2005.

\bibitem{Yazdani2018}
Kasra Yazdani, Austin Jones, Kevin Leahy, and Matthew~T. Hale.
\newblock {Differentially Private LQ Control}.
\newblock {\em arXiv preprint arXiv:1807.05082}, 2018.

\end{thebibliography}

\end{document}